\documentclass{amsart}
\usepackage{amsmath}
\usepackage{amsfonts}

\newtheorem{theorem}{Theorem}
\theoremstyle{plain}

\newtheorem{corollary}{Corollary}

\newtheorem{lemma}{Lemma}

\numberwithin{equation}{section}
\input{tcilatex}

\begin{document}
\title[\textbf{Frobenius-Euler numbers and polynomials}]{\textbf{The
Frobenius-Euler function and its applications}}
\author[\textbf{S. Araci}]{\textbf{Serkan Araci}}
\address{\textbf{University of Gaziantep, Faculty of Science and Arts,
Department of Mathematics, 27310 Gaziantep, TURKEY}}
\email{\textbf{mtsrkn@hotmail.com}}
\author[\textbf{D. Gao}]{\textbf{Deyao Gao}}
\address{\textbf{Yuren Lab., No. 8 Tongsheng Road, Changsha, P. R. China}}
\email{\textbf{13607433711@163.com}}
\author[\textbf{M. Acikgoz}]{\textbf{Mehmet Acikgoz}}
\address{\textbf{University of Gaziantep, Faculty of Science and Arts,
Department of Mathematics, 27310 Gaziantep, TURKEY}}
\email{\textbf{acikgoz@gantep.edu.tr}}

\begin{abstract}
In the present paper, we deal with Fourier-transformation of Frobenius-Euler
polynomials. We shall give its applications by using infinite series. Our
applications possess interesting properties which we state in this paper.

\vspace{2mm}\noindent \textsc{2010 Mathematics Subject Classification.}
11S80, 11B68.

\vspace{2mm}

\noindent \textsc{Keywords and phrases. }Frobenius-Euler numbers and
polynomials, Fourier transformation, infinite series.
\end{abstract}

\maketitle

\section{\textbf{Introduction}}

The ordinary Frobenius-Euler numbers are defined by means of the following
generating function:%
\begin{equation}
\sum_{n=0}^{\infty }H_{n}\left( u\right) \frac{t^{n}}{n!}=e^{H\left(
u\right) t}=\frac{1-u}{e^{t}-u}\text{.}  \label{equation 19}
\end{equation}%
where, in the umbral calculus, $H^{n}\left( u\right) $ is symbolically
replaced by $H_{n}\left( u\right) $ in the formal series expansion of 
\begin{equation*}
e^{tH\left( u\right) }=\sum_{n=0}^{\infty }H_{n}\left( u\right) \frac{t^{n}}{%
n!}\text{.}
\end{equation*}

From expression of this definition, we state the following%
\begin{equation}
\left( H\left( u\right) +1\right) ^{n}-uH_{n}\left( u\right) =\left\{ 
\begin{array}{cc}
1-u & \text{if }n=0, \\ 
0 & \text{if }n\in 
%TCIMACRO{\U{2115} }%
%BeginExpansion
\mathbb{N}
%EndExpansion
,%
\end{array}%
\right.   \label{equation 25}
\end{equation}%
where $%
%TCIMACRO{\U{2115} }%
%BeginExpansion
\mathbb{N}
%EndExpansion
$ denotes the set of positive integers.

From (\ref{equation 25}), we note that 
\begin{equation*}
H_{0}(u)=1,H_{1}\left( u\right) =-\frac{1}{1-u},H_{2}\left( u\right) =\frac{%
1+u}{\left( 1-u\right) ^{2}},\cdots .
\end{equation*}

The Frobenius Euler polynomials are also introduced as%
\begin{equation}
e^{xt}\frac{1-u}{e^{t}-u}=\sum_{n=0}^{\infty }\frac{t^{n}}{n!}H_{n}\left(
x,u\right) \text{.}  \label{equation 20}
\end{equation}

By (\ref{equation 19}) and (\ref{equation 20}), we can find the following%
\begin{equation*}
H_{n}\left( x,u\right) =\sum_{l=0}^{n}\binom{n}{l}x^{n-l}H_{l}\left(
u\right) \text{.}
\end{equation*}

By expression of (\ref{equation 25}), it is not difficult to show that the
recurrence relation for the Frobenius-Euler numbers as follows:%
\begin{equation}
\sum_{l=0}^{n}\binom{n}{l}H_{l}\left( u\right) -uH_{n}\left( u\right)
=\left( 1-u\right) \delta _{0,n}  \label{equation 26}
\end{equation}%
where $\delta _{m,n}$ is the Kronecker delta, is defined by%
\begin{equation*}
\delta _{m,n}=\left\{ 
\begin{array}{cc}
\text{ }1, & \text{if }m=n \\ 
0, & \text{ if }m\neq n.%
\end{array}%
\right. 
\end{equation*}

Thus, we easily procure the following: 
\begin{equation}
H_{n}\left( 1,u\right) -uH_{n}\left( u\right) =0\text{ }\left( \text{for }%
n\in 
%TCIMACRO{\U{2115} }%
%BeginExpansion
\mathbb{N}
%EndExpansion
\right) \text{.}  \label{equation 27}
\end{equation}

Thus, we arrive the following lemma.

\begin{lemma}
For $n\in 
%TCIMACRO{\U{2115} }%
%BeginExpansion
\mathbb{N}
%EndExpansion
$, we have%
\begin{equation*}
H_{n}\left( 1,u\right) =uH_{n}\left( u\right) \text{.}
\end{equation*}
\end{lemma}

Substituting $u=-1$ in the above lemma, it leads to%
\begin{equation*}
H_{n}\left( -1,u\right) =E_{n}\left( 1\right) =-E_{n}
\end{equation*}%
where $E_{n}$ is called Euler numbers, as is well-known, Euler numbers are
defined by the following generating function:%
\begin{equation*}
\sum_{n=0}^{\infty }E_{n}\frac{t^{n}}{n!}=\frac{2}{e^{t}+1}
\end{equation*}%
(for more informations on this subjects, see[1-20]).

Recently, Fourier transformation of the special functions have been studied
by many mathematicians (cf., \cite{Bayad 1}, \cite{J. Choi}, \cite{kim 2}, 
\cite{Kim 17}, \cite{kim 19}, \cite{Luo}). In \cite{Luo}, Luo gave Fourier
expansions of Apostol-Bernoulli and Apostol-Euler polynomials and derived
some integral representations of Apostol-Bernoulli and Apostol-Euler
polynomials by using Fourier expansions. After, Bayad \cite{Bayad 1}
introduced as theoretical identities of the Fourier transformation of the
Apostol-Bernoulli, Apostol-Euler and Apostol-Genocchi polynomials. Next, T.
Kim also defined the Euler function which is Fourier transformation of Euler
polynomials. We easily see that Kim's method is different from Bayad and
Luo. Actually, Kim's paper [3, pp. 131-136] motivated us to write this
paper. Thus, we also give Fourier transformation of Frobenius-Euler function
by using Kim's method. In this paper, we also show that this function is
related to Lerch trancendent $\Phi \left( z,s,a\right) $.

\section{\textbf{On the Frobenius-Euler function}}

In this section, we consider Frobenius-Euler function by using infinite
series. For $m\in 
%TCIMACRO{\U{2115} }%
%BeginExpansion
\mathbb{N}
%EndExpansion
$, the Fourier transformation of Frobenius-Euler function is introduced as%
\begin{equation}
H_{m}\left( x,u\right) =\sum_{n=-\infty }^{\infty }a_{n}^{\left( m\right)
}\left( u\right) e^{\left( 2n+1\right) \pi ix}\text{, }\left( a_{n}^{\left(
m\right) }\left( u\right) \in 
%TCIMACRO{\U{2102} }%
%BeginExpansion
\mathbb{C}
%EndExpansion
\right)   \label{equation 18}
\end{equation}%
where $%
%TCIMACRO{\U{2102} }%
%BeginExpansion
\mathbb{C}
%EndExpansion
$ denotes the set of complex numbers and%
\begin{equation}
a_{n}^{\left( m\right) }\left( u\right) =\int_{0}^{1}H_{m}\left( x,u\right)
e^{-\left( 2n+1\right) ix\pi }dx\text{.}  \label{equation 17}
\end{equation}

By applying some technical method on (\ref{equation 17}), we procure the
following%
\begin{eqnarray*}
a_{n}^{\left( m\right) }\left( u\right)  &=&\left[ \frac{H_{m+1}\left(
x,u\right) }{m+1}e^{-\left( 2n+1\right) ix\pi }\right] _{0}^{1}-\frac{\left(
2n+1\right) \pi i}{m+1}\int_{0}^{1}H_{m+1}\left( x,u\right) e^{-\left(
2n+1\right) ix\pi }dx \\
&=&-\frac{u+1}{m+1}H_{m+1}\left( u\right) -\frac{\left( 2n+1\right) \pi i}{%
m+1}a_{n}^{\left( m+1\right) }\left( u\right) \text{.}
\end{eqnarray*}

So from above, it leads to the following%
\begin{equation*}
a_{n}^{\left( m+1\right) }\left( u\right) =\left[ a_{n}^{\left( m\right)
}+\left( u+1\right) \frac{H_{m+1}\left( u\right) }{m+1}\right] \frac{m+1}{%
\left( \left( 2n+1\right) \pi i\right) }\text{.}
\end{equation*}

By continuing this process, becomes as follows: 
\begin{subequations}
\begin{gather}
a_{n}^{\left( m\right) }\left( u\right) =\left[ a_{n}^{\left( 1\right)
}\left( u\right) +\left( u+1\right) \frac{H_{2}\left( u\right) }{2}\right] 
\frac{m!}{\left( \left( 2n+1\right) \pi i\right) ^{m-1}}  \label{equation 16}
\\
+\left( u+1\right) \left[ \frac{1}{\left( 2n+1\right) \pi i}H_{m}\left(
u\right) +\frac{m}{\left( \left( 2n+1\right) \pi i\right) ^{2}}H_{m-1}\left(
u\right) +...+\frac{m!}{4!\left( \left( 2n+1\right) \pi i\right) ^{m-3}}%
\right] \text{.}  \notag
\end{gather}

We want to note that 
\end{subequations}
\begin{equation*}
\lim_{u\rightarrow -1}a_{n}^{\left( m\right) }\left( u\right) =a_{n}^{\left(
m\right) }\left( -1\right) :=a_{n}^{\left( m\right) }
\end{equation*}%
where $a_{n}^{\left( m\right) }\in 
%TCIMACRO{\U{2102} }%
%BeginExpansion
\mathbb{C}
%EndExpansion
$ is defined by Kim in \cite{Kim 1} as follows:%
\begin{equation*}
a_{n}^{\left( m\right) }=\frac{m!}{\left( \left( 2n+1\right) \pi i\right)
^{m-1}}a_{n}^{\left( 1\right) }\left( u\right) \text{.}
\end{equation*}

By using (\ref{equation 18}) and (\ref{equation 16}), we readily derive the
following%
\begin{equation*}
H_{m}\left( x,u\right) =\sum_{n=-\infty }^{\infty }\left\{ 
\begin{array}{c}
\left[ a_{n}^{\left( 1\right) }\left( u\right) +\left( u+1\right) \frac{%
H_{2}\left( u\right) }{2}\right] \frac{m!}{\left( \left( 2n+1\right) \pi
i\right) ^{m-1}}+\left( u+1\right)  \\ 
\times \left[ \frac{1}{\left( 2n+1\right) \pi i}H_{m}\left( u\right) +\frac{m%
}{\left( \left( 2n+1\right) \pi i\right) ^{2}}H_{m-1}\left( u\right) +...+%
\frac{m!}{4!\left( \left( 2n+1\right) \pi i\right) ^{m-3}}\right] 
\end{array}%
\right\} e^{\left( 2n+1\right) \pi ix}\text{.}
\end{equation*}

From this, we can state the following 
\begin{subequations}
\begin{gather*}
H_{m}\left( x,u\right) =\sum_{n=-\infty }^{\infty }\left[ a_{n}^{\left(
1\right) }\left( u\right) +\left( u+1\right) \frac{H_{2}\left( u\right) }{2}%
\right] \frac{m!e^{\left( 2n+1\right) \pi ix}}{\left( \left( 2n+1\right) \pi
i\right) ^{m-1}} \\
+\left( u+1\right) \sum_{n=-\infty }^{\infty }\left[ \frac{1}{\left(
2n+1\right) \pi i}H_{m}\left( u\right) +\frac{m}{\left( \left( 2n+1\right)
\pi i\right) ^{2}}H_{m-1}\left( u\right) +...+\frac{m!}{4!\left( \left(
2n+1\right) \pi i\right) ^{m-3}}H_{4}\left( u\right) \right] e^{\left(
2n+1\right) \pi ix}
\end{gather*}

After some calculations on the above equation, we have the following 
\end{subequations}
\begin{eqnarray*}
H_{m}\left( x,u\right)  &=&\sum_{n=-\infty }^{\infty }\left[ \frac{u+1}{u-1}%
\frac{1}{\left( 2n+1\right) \pi i}+\frac{2}{\left( \left( 2n+1\right) \pi
i\right) ^{2}}+\frac{1}{2}\left( \frac{u+1}{1-u}\right) ^{2}\right] \frac{%
m!e^{\left( 2n+1\right) \pi ix}}{\left( \left( 2n+1\right) \pi i\right)
^{m-1}} \\
&&+\left( u+1\right) \sum_{n=-\infty }^{\infty }\left[ \sum_{k=0}^{m-4}\frac{%
1}{\left( \left( 2n+1\right) \pi i\right) ^{k+1}}\frac{d^{k}}{dx^{k}}%
H_{m}\left( x,u\right) \mid _{x=0}\right] e^{\left( 2n+1\right) \pi ix}
\end{eqnarray*}

As a result, we conclude the following theorem.

\begin{theorem}
For $m\in 
%TCIMACRO{\U{2115} }%
%BeginExpansion
\mathbb{N}
%EndExpansion
$ and $0\leq x<1$, we have%
\begin{gather*}
H_{m}\left( x,u\right) =m!\sum_{n=-\infty }^{\infty }\left[ \frac{u+1}{u-1}%
\frac{1}{\left( 2n+1\right) \pi i}+\frac{2}{\left( \left( 2n+1\right) \pi
i\right) ^{2}}+\frac{1}{2}\left( \frac{u+1}{1-u}\right) ^{2}\right] \frac{%
e^{\left( 2n+1\right) \pi ix}}{\left( \left( 2n+1\right) \pi i\right) ^{m-1}}
\\
+\left( u+1\right) \sum_{k=0}^{m-4}\frac{1}{\left( 2\pi i\right) ^{k+1}}%
\frac{d^{k}}{dx^{k}}H_{m}\left( x,u\right) \mid _{x=0}\sum_{n=-\infty
}^{\infty }\frac{e^{\left( 2n+1\right) \pi ix}}{\left( n+\frac{1}{2}\right)
^{k+1}}\text{.}
\end{gather*}
\end{theorem}

Considering generating functions of Euler and Frobenius-Euler polynomials,
we reach the following corollary.

\begin{corollary}
Taking $u=-1$, we have Fourier transformation of Euler function, which is
defined by Kim in \cite{Kim 1} as follows:%
\begin{equation*}
H_{m}\left( x,-1\right) =E_{m}\left( x\right) =2m!\sum_{n=-\infty }^{\infty }%
\frac{e^{\left( 2n+1\right) \pi ix}}{\left( \left( 2n+1\right) \pi i\right)
^{m+1}}\text{ }\left( \text{for }0\leq x<1\right) \text{.}
\end{equation*}
\end{corollary}

The Lerch trancendent $\Phi \left( z,s,a\right) $ is the analytic
continuation of the series 
\begin{equation}
\Phi \left( z,s,a\right) =\sum_{n=0}^{\infty }\frac{z^{n}}{\left( n+a\right)
^{s}}  \label{equation 15}
\end{equation}%
which converges for $a\in 
%TCIMACRO{\U{2102} }%
%BeginExpansion
\mathbb{C}
%EndExpansion
\backslash 
%TCIMACRO{\U{2124} }%
%BeginExpansion
\mathbb{Z}
%EndExpansion
_{0}^{-}$, $s\in 
%TCIMACRO{\U{2102} }%
%BeginExpansion
\mathbb{C}
%EndExpansion
$ when $\left\vert z\right\vert <1$; $\Re \left( s\right) >1$ when $%
\left\vert z\right\vert =1$ where $%
%TCIMACRO{\U{2124} }%
%BeginExpansion
\mathbb{Z}
%EndExpansion
_{0}^{-}=%
%TCIMACRO{\U{2124} }%
%BeginExpansion
\mathbb{Z}
%EndExpansion
^{-}\cup \left\{ 0\right\} $, $%
%TCIMACRO{\U{2124} }%
%BeginExpansion
\mathbb{Z}
%EndExpansion
^{-}=\left\{ -1,-2,-3,...\right\} $. Lerch trancendent $\Phi \left(
z,s,a\right) $ is the proportional not only Riemann zeta funtion, Hurwitz
zeta function, the Dirichlet's eta function but also Dirichlet beta
function, the Legendre chi function, the polylogarithm, the Lerch zeta
function (for details, see \cite{Srivastava}, \cite{Acikgoz4}).

We now want to indicate that $\sum_{n=-\infty }^{\infty }\frac{e^{\left(
2n+1\right) \pi ix}}{\left( \left( 2n+1\right) \pi i\right) ^{m}}$ is
closely related to Lerch trancendent $\Phi \left( z,s,a\right) $. \ So, we
compute as follows:%
\begin{eqnarray*}
\sum_{n=-\infty }^{\infty }\frac{e^{\left( 2n+1\right) \pi ix}}{\left(
\left( 2n+1\right) \pi i\right) ^{m}} &=&\frac{1}{\left( 2\pi i\right) ^{m}}%
\sum_{n=-\infty }^{\infty }\frac{e^{\left( 2n+1\right) \pi ix}}{\left( n+%
\frac{1}{2}\right) ^{m}} \\
&=&\frac{1}{\left( 2\pi i\right) ^{m}}\sum_{n=-\infty }^{-1}\frac{e^{\left(
2n+1\right) \pi ix}}{\left( n+\frac{1}{2}\right) ^{m}}+\frac{1}{\left( 2\pi
i\right) ^{m}}\sum_{n=0}^{\infty }\frac{e^{\left( 2n+1\right) \pi ix}}{%
\left( n+\frac{1}{2}\right) ^{m}}
\end{eqnarray*}

After some applications on the above equation, we procure the following%
\begin{equation}
\sum_{n=-\infty }^{\infty }\frac{e^{\left( 2n+1\right) \pi ix}}{\left(
\left( 2n+1\right) \pi i\right) ^{m}}=-\frac{e^{\pi ix}}{\left( \pi i\right)
^{m}}+\frac{\left( -1\right) ^{m}e^{\pi ix}}{\left( 2\pi i\right) ^{m}}\Phi
\left( e^{-2\pi ix},m,-\frac{1}{2}\right) \text{.}  \label{equation 14}
\end{equation}

By using Theorem 1 and (\ref{equation 14}), we give the following theorem.

\begin{theorem}
The following equality holds true:%
\begin{gather*}
H_{m}\left( x,u\right) =m!\frac{u+1}{u-1}\left( -\frac{e^{\pi ix}}{\left(
\pi i\right) ^{m}}+\frac{\left( -1\right) ^{m}e^{\pi ix}}{\left( 2\pi
i\right) ^{m}}+\Phi \left( e^{-2\pi ix},m,-\frac{1}{2}\right) \right) \\
+2m!\left( -\frac{e^{\pi ix}}{\left( \pi i\right) ^{m+1}}+\frac{\left(
-1\right) ^{m+1}e^{\pi ix}}{\left( 2\pi i\right) ^{m+1}}+\Phi \left(
e^{-2\pi ix},m+1,-\frac{1}{2}\right) \right) \\
+\frac{1}{2}\left( \frac{u+1}{u-1}\right) ^{2}\left( -\frac{e^{\pi ix}}{%
\left( \pi i\right) ^{m-1}}+\frac{\left( -1\right) ^{m-1}e^{\pi ix}}{\left(
2\pi i\right) ^{m-1}}+\Phi \left( e^{-2\pi ix},m-1,-\frac{1}{2}\right)
\right) \\
+\left( u+1\right) \sum_{k=0}^{m-4}\frac{d^{k}}{dx^{k}}H_{m}\left(
x,u\right) \mid _{x=0}\left( -\frac{e^{\pi ix}}{\left( \pi i\right) ^{k+1}}+%
\frac{\left( -1\right) ^{k+1}e^{\pi ix}}{\left( 2\pi i\right) ^{k+1}}+\Phi
\left( e^{-2\pi ix},k+1,-\frac{1}{2}\right) \right) \text{.}
\end{gather*}
\end{theorem}

For $u=-1$ on the above theorem, we have the following corollary.

\begin{corollary}
The following identity%
\begin{equation*}
E_{m}\left( x\right) =2m!\left( -\frac{e^{\pi ix}}{\left( \pi i\right) ^{m+1}%
}+\frac{\left( -1\right) ^{m+1}e^{\pi ix}}{\left( 2\pi i\right) ^{m+1}}+\Phi
\left( e^{-2\pi ix},m+1,-\frac{1}{2}\right) \right)
\end{equation*}%
is true.
\end{corollary}

Setting $x=1$ in Theorem 1, we obtain%
\begin{gather}
H_{m}\left( 1,u\right) =-m!\sum_{n=-\infty }^{\infty }\left[ \frac{u+1}{u-1}%
\frac{1}{\left( 2n+1\right) \pi i}+\frac{2}{\left( \left( 2n+1\right) \pi
i\right) ^{2}}+\frac{1}{2}\left( \frac{u+1}{1-u}\right) ^{2}\right] \frac{1}{%
\left( \left( 2n+1\right) \pi i\right) ^{m-1}}  \label{equation 13} \\
-\left( u+1\right) \sum_{k=0}^{m-4}\frac{1}{\left( 2\pi i\right) ^{k+1}}%
\frac{d^{k}}{dx^{k}}H_{m}\left( x,u\right) \mid _{x=0}\sum_{n=-\infty
}^{\infty }\frac{1}{\left( n+\frac{1}{2}\right) ^{k+1}}\text{.}  \notag
\end{gather}

By expressions of (\ref{equation 13}) and Lemma 1, we easily see the
following corollary.

\begin{corollary}
The following identity holds true:%
\begin{gather*}
uH_{m}\left( u\right) =-m!\sum_{n=-\infty }^{\infty }\left[ \frac{u+1}{u-1}%
\frac{1}{\left( 2n+1\right) \pi i}+\frac{2}{\left( \left( 2n+1\right) \pi
i\right) ^{2}}+\frac{1}{2}\left( \frac{u+1}{1-u}\right) ^{2}\right] \frac{1}{%
\left( \left( 2n+1\right) \pi i\right) ^{m-1}} \\
-\left( u+1\right) \sum_{k=0}^{m-4}\frac{1}{\left( 2\pi i\right) ^{k+1}}%
\frac{d^{k}}{dx^{k}}H_{m}\left( x,u\right) \mid _{x=0}\sum_{n=-\infty
}^{\infty }\frac{1}{\left( n+\frac{1}{2}\right) ^{k+1}}\text{.}
\end{gather*}
\end{corollary}

Now, by using Kim's method in \cite{Kim 1}, we discover the following 
\begin{subequations}
\begin{gather}
\frac{1}{1-ue^{-t}}=\sum_{n=0}^{\infty }u^{n}e^{-nt}=\sum_{n=0}^{\infty
}u^{n}\left( e^{-t}\right) ^{n}=\sum_{n=0}^{\infty }u^{n}\left(
\sum_{k=0}^{\infty }\left( -1\right) ^{k}\frac{t^{k}}{k!}\right) ^{n}
\label{equation 12} \\
=\sum_{n=0}^{\infty }\left( \sum_{a_{1}+a_{2}+...+a_{n}=n}\frac{n!}{\left(
a_{1}\right) !\left( a_{2}\right) !...}\frac{\left( -1\right)
^{a_{1}+2a_{2}+...}}{\left( 1!\right) ^{a_{1}}\left( 2!\right) ^{a_{2}}...}%
\right) t^{a_{1}+2a_{2}+...}  \notag
\end{gather}

Let $p\left( i,j\right) :a_{1}+2a_{2}+...=i,a_{1}+a_{2}+...=j$, from
expression of (\ref{equation 12}), we compute as follows: 
\end{subequations}
\begin{subequations}
\begin{gather}
\frac{1}{1-ue^{-t}}=\sum_{m=0}^{\infty }\left(
\sum_{n=0}^{m}u^{n}\sum_{p\left( m,n\right) }\frac{n!}{\left( a_{1}\right)
!\left( a_{2}\right) !...\left( a_{m}\right) !}\frac{\left( -1\right)
^{a_{1}+2a_{2}+...+ma_{m}}}{\left( 1!\right) ^{a_{1}}\left( 2!\right)
^{a_{2}}...\left( m!\right) ^{a_{m}}}\right) t^{a_{1}+2a_{2}+...+ma_{m}} 
\notag \\
=\sum_{m=0}^{\infty }\left( -1\right) ^{m}\left(
\sum_{n=0}^{m}n!u^{n}\sum_{p\left( m,n\right) }\frac{m!}{\left( a_{1}\right)
!\left( a_{2}\right) !...\left( a_{m}\right) !}\frac{\left( -1\right) ^{m}}{%
\left( 1!\right) ^{a_{1}}\left( 2!\right) ^{a_{2}}...\left( m!\right)
^{a_{m}}}\right) \frac{t^{m}}{m!}  \notag \\
=\sum_{m=0}^{\infty }\left[ \left( -1\right)
^{m}\sum_{n=0}^{m}n!u^{n}s_{2}\left( m,n\right) \right] \frac{t^{m}}{m!}%
\text{.}  \label{equation 11}
\end{gather}%
where $s_{2}\left( m,n\right) $ is the second kind stirling number.

Via the definition of Frobenius-Euler numbers, we readily derive the
following 
\end{subequations}
\begin{eqnarray}
\frac{1}{1-ue^{-t}} &=&\frac{u^{-1}}{u^{-1}-1}\frac{1-u^{-1}}{e^{-t}-u^{-1}}
\label{equation 10} \\
&=&\frac{1}{1-u}\sum_{m=0}^{\infty }\left( -1\right) ^{m}H_{m}\left(
u^{-1}\right) \frac{t^{m}}{m!}\text{.}  \notag
\end{eqnarray}

By comparing the coefficients of $\frac{t^{n}}{n!}$ on the both sides of (%
\ref{equation 11}) and (\ref{equation 10}), we reach the following theorem.

\begin{theorem}
\bigskip The following equality holds true:%
\begin{equation*}
\frac{1}{1-u}H_{m}\left( u^{-1}\right) =\sum_{n=0}^{m}n!u^{n}s_{2}\left(
m,n\right) \text{.}
\end{equation*}
\end{theorem}

\begin{corollary}
Substituting $u=-1$ in the above theorem, we get the following, which is
defined by Kim \cite{Kim 1}%
\begin{equation*}
E_{m}=2\sum_{n=0}^{m}n!\left( -1\right) ^{n}s_{2}\left( m,n\right) \text{.}
\end{equation*}
\end{corollary}

\end{document}